\documentstyle[12pt,euscript]{amsart}
\def\mysavedown#1{\edef\mysubs{\mysubs#1}}
\def\mysaveup#1{\edef\mysups{\mysups#1}}
\def\mydown#1{{\mytensor}_{\vphantom{\mysubs}#1}}
\def\myup#1{{\mytensor}^{\vphantom{\mysups}#1}}
\def\tensor#1#2{
  #1
  \def\mytensor{\vphantom{#1}}
  \def\mysubs{\relax}
  \def\mysups{\relax}
  \let\down=\mysavedown
  \let\up=\mysaveup
  #2
  \let\down=\mydown
  \let\up=\myup
  #2
  }
\renewcommand{\phi}{\varphi}
\renewcommand{\epsilon}{\varepsilon}

\newcommand{\C}{\bold C}
\newcommand{\R}{\bold R}
\newcommand{\cross}{\mathbin{\times}}
\newcommand{\vhat}{{\widehat V}}

\newcommand{\shat}{{\widehat{\scr S}}}
\newcommand{\Phihat}{{\widehat\Phi}}
\newcommand{\Psihat}{{\widehat\Psi}}

\newtheorem{theorem}{Theorem}[section]
\newtheorem{lemma}[theorem]{Lemma}

\newtheorem{conjecture}[theorem]{Conjecture}
\newtheorem{bigtheorem}{Theorem}[section]

\let\scr=\EuScript

\def\crn#1#2{{\vcenter{\vbox{
        \hbox{\kern#2pt \vrule width.#2pt height#1pt
           }
          \hrule height.#2pt}}}}
\def\intprod{\mathchoice\crn54\crn54\crn{3.75}3\crn{2.5}2}
\def\into{\mathbin{\intprod}}
\def\tprod{\mathbin{\otimes}}

\begin{document}
\title
{CR manifolds with noncompact connected automorphism groups}
\author{John M. Lee}
\thanks{Research supported in part by National
Science Foundation grant DMS 91-01832.}
\address{%
Department of Mathematics, GN-50\\
  University of Washington \\
  Seattle, WA 98195 }
\email{lee@@math.washington.edu}
\maketitle
\begin{abstract}
The main result of this paper is that the identity component of the
automorphism group of a compact, connected, strictly pseudoconvex CR
manifold is compact unless the manifold is CR equivalent to the standard
sphere.  In dimensions greater than 3, it has been pointed out by D. Burns
that this result follows from known results on biholomorphism groups of
complex manifolds with boundary and the fact that any such CR manifold $M$
can be realized as the boundary of an analytic variety.  When $M$ is
3-dimensional, Burns's proof breaks down because abstract CR 3-manifolds
are generically not realizable as boundaries.  This paper provides an
intrinsic proof of compactness that works in any dimension.
\end{abstract}

\section {Introduction}

The purpose of this paper is to prove the following theorem.

\begin{bigtheorem}\label{thma}
Let $M$ be a compact, connected, strictly pseudoconvex CR manifold of
dimension $2n+1\ge 3$.  Then the identity component $\scr A_0(M)$ of the
group $\scr A(M)$ of CR automorphisms of $M$ is compact unless $M$ is
globally CR equivalent to the $(2n+1)$-sphere with its standard CR
structure.
\end{bigtheorem}

In dimensions greater than 3, this result follows from known results on
biholomorphism groups of complex manifolds with boundary.  In fact, a
stronger result is true in that case.  The following is due to Dan
Burns, although it has never been published:

\begin{theorem}\label{burnsthm}
\rom{(}D. Burns\rom{)} 
Suppose $M$ is a compact, connected, strictly pseudoconvex CR manifold of
dimension $2n+1\ge 5$.  The full CR automorphism group $\scr A(M)$ is
compact unless $M$ is globally CR equivalent to $S^{2n+1}$ with its
standard CR structure.
\end{theorem}

The idea is that any such $M$ can be realized as the boundary of an
analytic variety whose biholomorphism group is isomorphic to the group of
CR automorphisms of $M$.  It follows from results of B.  Wong \cite{Wong},
J.-P. Rosay \cite{Rosay}, and D. Burns and S. Shnider \cite{BS} (with a
slight extra argument to deal with singular varieties) that any such
variety with noncompact biholomorphism group is biholomorphic to the unit
ball, so $M$ must be the sphere.  Theorem \ref{burnsthm} does not appear
explicitly in the literature, so for completeness we give a sketch of
Burns's proof in Section 2.

When $M$ is 3-dimensional, this proof breaks down because abstract CR
3-manifolds are generically not realizable as boundaries.  The best
previously known result is the following theorem of S. Webster.  We say
that a CR manifold is {\it locally spherical} if it is locally CR
equivalent to the sphere with its standard CR structure.

\begin{theorem}\label{websterthm}
\rom{(}S. Webster \cite{Webster}\rom{)} 
If $M$ is a compact, connected, strictly pseudoconvex CR manifold, and
$\scr A_0(M)$ is noncompact, then $M$ is locally spherical.
\end{theorem}

To put these results into perspective, is it useful to consider conformal
Riemannian geometry as a guide, since there is a strong analogy between
conformal and CR geometry.  The analogous result for conformal manifolds is
the following theorem.  Recall that a Lie group $G$ acts {\it properly} on
a space $X$ if the map $G\cross X\to X\cross X$ given by $(g,x)\mapsto
(g\cdot x, x)$ is proper.

\begin{theorem} \label{confthm} 
Let $X$ denote a Riemannian manifold, $\scr C(X)$
its group of conformal diffeomorphisms, and $\scr C_0(X)$ the identity
component of $\scr C(X)$.
\begin{enumerate}
\item\label{obatathm}
\rom{(}M. Obata \cite{Obata}, U. Pinkall, and J. Lafontaine
\cite{Lafontaine}\rom{)} 
If $X$ is compact, then $\scr C_0(X)$ is compact unless $X$ is conformally
equivalent to the sphere with its standard metric.
\item
\rom{(}J. Ferrand \cite{F1}\rom{)} 
The same is true with $\scr C_0(X)$ replaced by $\scr C(X)$.
\item\label{noncompact}
\rom{(}J. Ferrand \cite{F2}\rom{)} 
If $X$ is non-compact, then $\scr C(X)$ acts properly unless $X$ is
conformally equivalent to $\R^n$ with the Euclidean metric.
\end{enumerate}
\end{theorem}

These results have an interesting history.  Theorem
\ref{confthm}(\ref{obatathm}) was originally claimed by Obata in
\cite{Obata}, and is commonly attributed to him.  However, that
paper contained a gap, which was repaired in
\cite{Lafontaine}, based on an argument due to Pinkall.
Theorem \ref{confthm}(\ref{noncompact}) was claimed by V. Alekseevskii
\cite{A1,A2}.  However, the proof in \cite{A1} 
apparently contains a serious error (Theorem 4 is false), so the result was
in question until the appearance of \cite{F2}.  The expository paper
\cite{Gutschera} by R. Gutschera gives an excellent survey of results in
the conformal category.

By analogy with the conformal case, it is reasonable to make the following
conjecture.

\begin{conjecture}\label{conj}
If $M$ is a connected, strictly pseudoconvex CR manifold, $\scr A(M)$ acts
properly unless $M$ is CR equivalent to the sphere or the Heisenberg group
with its standard CR structure.
\end{conjecture}

In case $M$ is compact, properness of the action implies compactness of
$\scr A(M)$, so this conjecture includes as a special case the conjecture
that $\scr A(M)$ is compact when $M$ is a compact strictly pseudoconvex
3-manifold other than the sphere.  When $M$ is noncompact, one can show as
in \cite{Gutschera} that $\scr A(M)$ acts properly if and only if it
preserves a pseudohermitian structure (see \S\ref{defsection} for
definitions).

P. Pansu has recently pointed out \cite{Pansu} that Ferrand's methods can
be extended to the CR case to show that $\scr A(M)$ is precompact in the
$C^0$ topology.  But that approach is not yet strong enough to prove the
full strength of Theorem \ref{thma} or Theorem \ref{burnsthm}, because it
is not known in general whether a nonconstant uniform limit of CR
diffeomorphisms is a diffeomorphism.

The proof of Theorem \ref{thma} is carried out along the lines of the proof
of Theorem \ref{confthm}(\ref{obatathm}).  In fact, as Webster already
observed in \cite{Webster}, all but one step of that proof goes through in
the CR case with little difficulty.  Webster showed the following:

\begin{theorem}\label{websterfixptthm}
\rom{(}S. Webster \cite{Webster}\rom{)} 
If $M$ is a compact, connected, strictly pseudoconvex, locally spherical CR
manifold, and there exists a closed, noncompact one-parameter subgroup
$G_1\subset \scr A_0(M)$ with a fixed point, then $M$ is globally CR
equivalent to $S^{2n+1}$.
\end{theorem}
\noindent
(The hypothesis that $G_1$ is closed was omitted from the statement of the
theorem in \cite{Webster}, but it is clearly necessary.)  Thus the only
tricky part of Theorem \ref{thma} is proving that a closed, noncompact
1-parameter subgroup of $\scr A_0(M)$ has a fixed point.  The bulk of this
paper is devoted to proving this fixed-point result, which we state as a
separate theorem.

\begin{bigtheorem}\label{fixptthm}
Let $M$ be a compact, connected, strictly pseudoconvex CR manifold of
dimension $2n+1 \ge 3$, and $G_1\subset \scr A_0(M)$ a closed, noncompact
one-parameter subgroup.  Then $G_1$ has a fixed point.
\end{bigtheorem}

The corresponding step in the conformal case \cite{Obata} was trivial: if
the infinitesimal generator $X$ of $G_1$ never vanishes, then one can
rescale the metric so that $X$ has norm 1, and then $G_1$ preserves the
rescaled metric, so compactness follows easily.  In the CR case, the
analogous argument only allows us to conclude that $X$ is tangent to the
contact bundle $H$ somewhere.  To draw the stronger conclusion that $X$
vanishes somewhere, we must carefully analyze the set $S$ where $X$ is
tangent to $H$.  This turns out to be a smoothly embedded compact
hypersurface, which carries a Riemannian metric preserved by $G_1$.  This
implies that any sequence of elements of $G_1$ has a subsequence that
converges along $S$.  A further calculation shows that convergence along
$S$ entails convergence of the 2-jets, from which global convergence
follows.

Since there is no extra work involved, the proof is carried out in all
dimensions, thereby providing an independent (and considerably more
elementary) proof of Theorem \ref{burnsthm} in the special case when the
identity component $\scr A_0(M)$ is noncompact.

Here is the proof of Theorem \ref{thma}, given Theorem \ref{fixptthm}.
Assuming that $M$ is as in the statement of Theorem \ref{thma} and $\scr
A_0(M)$ is noncompact, Webster's Theorem \ref{websterthm} implies that $M$
is locally spherical.  It follows from the properties of the Cartan
connection constructed by S. S. Chern and J. K.  Moser \cite{CM} and
standard results on $G$-structures that $\scr A(M)$ has a unique smooth
manifold structure making it into a Lie transformation group.  (In fact,
the topology on $\scr A(M)$ can be taken to be that of $C^2$ convergence on
$M$.)  If $\scr A_0(M)$ is not compact, then by an old theorem of D.
Montgomery and L.  Zippin \cite{MZ} $\scr A_0(M)$ has a closed
one-parameter subgroup $G_1$ which is isomorphic to $\bold R$.  By Theorem
\ref{fixptthm}, $G_1$ has a fixed point.  But then Theorem
\ref{websterfixptthm} implies that $M$ is globally CR equivalent to
$S^{2n+1}$.  Thus the theorem is proved.

An important application of this compactness result is to the construction
of local slices for the action of the contact diffeomorphism group on the
set of CR structures on a 3-manifold.  This is carried out in a joint paper
with J.-H. Cheng \cite{CL}.

In Section \ref{defsection}, we introduce our notation and review some
facts from the theory of CR and pseudohermitian manifolds.  At the end of
the section, we sketch the proof of Theorem \ref{burnsthm}.  In Section
\ref{fixptsection}, we prove the fixed point theorem, Theorem \ref{fixptthm}.

I would like to thank all the people with whom I have had useful
discussions about this work, especially Jih-Hsin Cheng, Dan Burns, Robert
Gutschera, Lee Stout, and Lutz Bungart.

\section{Background}\label{defsection}

Throughout this paper, we use the notation and terminology of
\cite{pseudo} unless otherwise specified; we refer the reader there for
basic notions of CR geometry not explained here.  Suppose $M$ is a
hypersurface-type CR manifold of dimension $2n+1$.  This means $M$ is
endowed with a smooth $n$-dimensional complex subbundle $\scr H\subset
TM\tprod \C$ which satisfies $\scr H\cap
\scr H = \{0\}$ and which is formally integrable:
$[\Gamma(\scr H), \Gamma(\scr H)]\subset\Gamma(\scr H)$.  The real bundle
$H = \operatorname{Re}(\scr H\oplus \overline{\scr H})\subset TM$ carries a
complex structure map $J\colon H\to H$ satisfying $J^2=-1$, and $\scr H$ is
just the $i$-eigenspace of the complexification of $J$.

An orientation of $TM$ together with the orientation of $H$ induced by $J$
automatically induces an orientation on the annihilator of $H$ in $T^*M$,
so if $M$ is orientable there exists a global real 1-form $\theta$ whose
kernel at each point is $H$.  Once such a form $\theta$ is chosen, the {\it
Levi form} determined by $\theta$ is the symmetric bilinear form on $H$
defined by
\begin{displaymath}
{\left< V, W \right>}_\theta = d\theta(V,J W).
\end{displaymath}
The same formula, extended by complex bilinearity, gives a complex-bilinear
form on $H\tprod \bold C$ which is Hermitian on $\scr H\cross
\overline{\scr H}$.  A standard
computation shows that the Levi form changes conformally if $\theta$ is
changed, so its signature is a CR invariant of $M$.  If the Levi form is
positive definite, the CR structure is said to be {\it strictly
pseudoconvex}, and in that case $\theta$ is a contact form.  A strictly
pseudoconvex CR structure together with a given contact form is called a
{\it pseudohermitian structure}.

On a pseudohermitian manifold, the Levi form yields a norm on all (real or
complex) tensor bundles over $H$, denoted ${|\cdot|}_\theta$, and a {\it
characteristic vector field} $T$, defined by
\begin{equation}\label{defT}
T\into\theta=1,\qquad  T\into d\theta=0.
\end{equation}
There is also a natural linear connection, the {\it pseudohermitian
connection} \cite{Websterthesis,Websterconn,Tanaka}, and a natural
Riemannian metric $g_\theta$ called the {\it Webster metric}
\cite{Websterconn}, characterized as the unique inner product on $TM$ that
restricts to the Levi form on $H$, and for which $T$ is a unit vector
orthogonal to $H$.  If we extend $J$ to an endomorphism $J_\theta$ of $TM$
by declaring $J_\theta T = 0$, then $g_\theta$ can be written
\begin{displaymath}
g_\theta(X,Y) = d\theta(X,J_\theta Y) + \theta(X)\theta(Y).
\end{displaymath}

For local computations on a pseudohermitian manifold, it is useful to
choose a complex local frame $(T,Z_\alpha, Z_{\bar\alpha})$, where
$\{Z_\alpha: \alpha=1,\dots,n\}$ forms a basis for $\scr H$ at each point,
$Z_{\bar\alpha} = \overline {Z_\alpha}$, and $T$ is the characteristic
vector field of $\theta$.  We let $\{\theta, \theta^\alpha,
\theta^{\bar\alpha}\}$ denote the dual coframe.  The components of a tensor
with respect to this frame are denoted by subscript and/or superscript
indices using the summation convention, with a zero index referring to the
$T$ direction.  For example, if $\eta$ is a 1-form, we can write locally
\begin{displaymath}
\eta = \eta_0\theta + \eta_\alpha\theta^\alpha + \eta_{\bar\alpha}
\theta^{\bar\alpha}.
\end{displaymath}
The components of the Levi form in such a local frame are
$h_{\alpha\bar\beta}$, where
\begin{displaymath}
d\theta = ih_{\alpha\bar\beta}\theta^\alpha\wedge\theta^{\bar\beta}.
\end{displaymath}
We use the Hermitian matrix $h_{\alpha\bar\beta}$ and its inverse
$h^{\alpha\bar\beta}$ to lower and raise indices in the usual way.  The
pseudohermitian connection is determined by a matrix of complex-valued
1-forms $\omega_\beta{}^\alpha$ satisfying
\begin{equation}\label{dthetaa}
d\theta^\alpha = \theta^\beta\wedge \omega_\beta{}^\alpha +
A^\alpha{}_{\bar\beta}\theta\wedge\theta^{\bar\beta},
\end{equation}
where $A^\alpha{}_{\bar\beta}$ are the components of the Webster torsion
tensor.
The pseudohermitian curvature tensor is denoted by $\tensor R
{\down\alpha \up \rho \down{\beta\bar\gamma}}$, and the pseudohermitian
Ricci tensor by $R_{\alpha\bar\gamma} = \tensor R {\down\alpha \up \beta
\down{\beta\bar\gamma}}$.

Components of covariant derivatives are indicated with indices preceded by
a semicolon.  The most important commutation relations for covariant
derivatives are the following (cf.\ \cite{pseudo}): if $f$ is a smooth
function on $M$, then
\begin{gather}\label{comm2}
f_{;\alpha\bar\beta} - f_{;\bar\beta\alpha} =
ih_{\alpha\bar\beta}f_{;0};
\qquad f_{;0\alpha}-f_{;\alpha 0} =
A_{\alpha\beta}f^{;\beta};\\
\label{comm3}
f_{;\alpha\beta\bar\gamma} - f_{;\alpha\bar\gamma\beta} =
ih_{\beta\bar\gamma} f_{;\alpha 0} +
\tensor R {\down\alpha \up \rho \down{\beta\bar\gamma}} f_{;\rho}.
\end{gather}

For a smooth function $f$, $\overline\partial_b f$ is the restriction to
$\overline{\scr H}$ of $df$.  On a pseudohermitian manifold, we can
identify $\overline\partial_b f$ as an honest differential form by
stipulating that $T\into \overline\partial_b f=0$, so in terms of a local
frame we have
\begin{displaymath}
\overline\partial_b f = f_{;\bar\alpha}\theta^{\bar\alpha}.
\end{displaymath}
The norm of $\overline\partial_b f$ is $|\overline\partial_b f|_\theta^2 =
f_{;\bar\alpha} f^{;\bar\alpha}$.

A real vector field $X$ on $M$ whose flow acts by CR automorphisms is
called a {\it CR vector field}.  In particular, since the flow preserves
the contact bundle $H$, any such vector field is an infinitesimal contact
automorphism.  It is well known (cf.\ \cite{Gray,CL})
that any such vector field is of the form
\begin{equation}\label{contactvf}
X = H^\theta_f - fT
\end{equation}
for some real-valued smooth function $f$, where $H^\theta_f$ is the {\it
contact Hamiltonian field} of $f$, defined by $H^\theta_f\into \theta=0$,
$H^\theta_f\into d\theta = df - (Tf)\theta$.  In terms of a local frame,
\begin{equation}\label{hf}
H^\theta_f = i f^{;\bar\alpha}Z_{\bar\alpha} - i f^{;\alpha}Z_{\alpha}.
\end{equation}
If $X$ is a CR vector field,
then in addition $X$ satisfies
\begin{displaymath}
L_X\theta^{\bar\beta} \equiv 0 \pmod {\theta, \theta^{\bar\alpha}}
\end{displaymath}
which implies (using \ref{dthetaa})
\begin{displaymath}
0 = (L_X\theta^{\bar\beta})(Z_{\alpha})
= i f^{;\bar\beta}{}_{\alpha} - A^{\bar\beta}{}_{\alpha}f,
\end{displaymath}
or, lowering indices,
\begin{equation}\label{fab}
f_{;\beta\alpha} = - i A_{\beta\alpha}f.
\end{equation}

We conclude this section by sketching a proof of Theorem \ref{burnsthm}.
Let $M$ be as in the statement of the theorem.  By Boutet de Monvel's
embedding theorem \cite{BdM}, there exist smooth global CR-holomorphic
functions $(z^1,\dots, z^m)$ on $M$ that define a CR embedding of $M$ into
$\bold C^m$ for some  $m$.  Then by the ``filling-in'' theorem of R.
Harvey and B.  Lawson \cite{HL}, $M$ bounds a unique compact complex
analytic variety $V\subset\bold C^m$, which is smooth except possibly at a
finite set $\scr S$
of isolated singular points in its interior.  

Let $\pi\colon \vhat\to V$ denote the normalization of $V$ \cite{Gunning};
thus $\vhat$ is an abstract normal Stein analytic space with a finite
singular set $\shat\subset \pi^{-1}(\scr S)$, and $\pi$ is a biholomorphism
on the complement of the finite set $\pi^{-1}(\scr S)$.  It follows that
$\vhat$ also has $M$ as smooth boundary.

Let $\scr B(\vhat)$ denote the group of biholomorphisms of $\vhat$.
Equipped with the compact-open topology, $\scr B(\vhat)$ is a Lie group
\cite{Fujimoto}.  I claim that $\scr A(M)$ and $\scr B(\vhat)$
are topologically isomorphic.  Let $\phi\in\scr A(M)$ be any CR
automorphism of $M$.  Since $M=\partial V$ is strictly pseudoconvex, the
coordinate functions of $\phi\colon M\to M\subset \C^m$ extend uniquely to
a neighborhood of $M$ in $V$; then by Hartogs's theorem for Stein varieties
\cite[p.\ 228]{Gunningbook}, $\phi$ extends to a weakly holomorphic map
$\Phi\colon V\to V$.  By normality, the weakly holomorphic map
$\Phi\circ\pi\colon \vhat\to V$ extends holomorphically to all of $\vhat$.
Moreover, away from the proper subvariety $(\Phi\circ\pi)^{-1}(\scr S)$,
this map lifts to a map $\Phihat = \pi^{-1} \circ \Phi \circ \pi \colon
\vhat\to\vhat$, and by normality again, the lifted map $\Phihat$ extends to
all of $\vhat$.  Arguing similarly for $\phi^{-1}$, we obtain a map
$\Psihat\colon\vhat\to\vhat$ such that $\Psihat\circ\Phihat$ and
$\Phihat\circ\Psihat$ are the identity near the boundary of $\vhat$, and by
analytic continuation globally; thus $\Phihat$ is a biholomorphism.

This construction yields a group homomorphism $E\colon \scr A(M)\to\scr
B(\vhat)$.  Fefferman's extension theorem \cite{Fefferman} (suitably
localized near the boundary as in \cite{Forst}) shows that any
biholomorphism $\Phihat$ of $\vhat$ extends smoothly to $M$, so $E$ is
bijective.  Since a Stein variety is holomorphically convex, convergence of
$\{\phi_j\}\subset \scr A(M)$ implies uniform convergence of
$\{\Phihat_j\}$, so $E$ is continuous.  Since any continuous bijection
between Lie groups is a diffeomorphism, it follows that $\scr A(M)$ and
$\scr B(\vhat)$ are smoothly isomorphic.

Suppose $\scr A(M)$ is noncompact.  By the discussion above, $\scr
B(\vhat)$ is also noncompact.  If we knew $\vhat$ were smooth,
noncompactness of $\scr B(\vhat)$ would imply by \cite[Thm.\ II]{BS} (cf.\
also \cite{Wong,Rosay}) that $\vhat$ is biholomorphically equivalent to the
$(n+1)$-ball.  By Fefferman's theorem once again, the biholomorphism
extends to a CR equivalence between $M$ and the sphere, thus proving the
theorem.  Thus to complete the proof it suffices to show that $\shat$ is
empty when $\scr B(\vhat)$ is noncompact.

Suppose on the contrary that $\shat$ is not empty and $\scr B(\vhat)$ is
noncompact.  Then so is the subgroup $\scr B_0(\vhat)$ consisting of
biholomorphisms that fix $\shat$ pointwise, since $\scr B_0(\vhat)$ has
finite index in $\scr B(\vhat)$ due to the finiteness of $\shat$.  Let
$\{\Phihat_j\}\subset \scr B_0(\vhat)$ be a sequence with no subsequence
that converges in the compact-open topology to a biholomorphism.  By
applying Montel's theorem for varieties \cite{Gunning-Vitali,Bungart} to
the coordinate functions of $\pi\circ\Phihat_j$, one can show that a
subsequence converges uniformly on compact subsets to a holomorphic map
$\Phihat\colon\vhat\to \vhat\cup M$.  Passing to a smaller subsequence and
arguing similarly for $\Phihat_j^{-1}$, we may assume also that
$\Phihat_j^{-1}$ converges to a holomorphic map $\Psihat\colon\vhat\to
\vhat\cup M$.

If $\Phihat(p)\in M$ for some $p\in \vhat$, then letting $f$ be a
holomorphic peak function on $\vhat\cup M$ taking its maximum modulus at
$\Phi(p)$, $f\circ\Phi$ would take an interior maximum at $p$, so would be
constant by the maximum principle; but this contradicts the fact that each
$\Phihat_j$ fixes $\shat$.  Thus $\Phihat$ maps $\vhat$ to $\vhat$; a
similar comment applies to $\Psihat$.  Then uniform convergence of
$\Phihat_j$ and $\Phihat_j^{-1}$ on a compact neighborhood of $\shat$
implies $\Psihat\circ\Phihat$ and $\Phihat\circ\Psihat$ are the identity on
an open set, and by analytic continuation globally, so $\Phihat_j$
converges to a biholomorphism, which is a contradiction.  This completes
the proof.

\section{A fixed-point theorem}\label{fixptsection}

In this section we prove Theorem \ref{fixptthm}, from which Theorem
\ref{thma} follows as described in the Introduction.

We need to show that the infinitesimal generator $X$ of $G_1$ has a zero.
The proof will be by contradiction.  We begin by noting that we can assume
without loss of generality that $M$ is orientable: if not, just lift the CR
structure and $X$ to the orientable double cover, and the lifted vector
field will have zeros if and only if $X$ does.  When $M$ is orientable, it
is possible to choose a global contact form $\theta$.

Before proving Theorem \ref{fixptthm}, we establish some preliminary
lemmas.  Throughout this section, $M$ will be as in the statement of the
theorem and $X$ will be the infinitesimal generator of $G_1$, with the
additional assumptions that $M$ is orientable and $X$ has no zeros on $M$.

\begin{lemma}
There is a nonempty, compact, embedded hypersurface $S\subset M$ along
which $X$ is tangent to both $H$ and $S$.
\end{lemma}

\begin{pf}
Choose a global contact form $\theta$ on $M$.  Since $X$ is an
infinitesimal contact diffeomorphism, we can write $X$ in the form
(\ref{contactvf}), where $f=-\theta(X)$ satisfies (\ref{fab}).  If $f$
vanishes nowhere, then replacing $\theta$ by $\widetilde\theta =
(1/f)\theta$, we see that $\widetilde f = -\widetilde\theta(X) = 1$.  This
implies $H^{\tilde\theta}_{\tilde f} = 0$, $X = -\widetilde T$, and
therefore
\begin{displaymath}
L_X\widetilde\theta = -d(\widetilde T\into\widetilde\theta) -
\widetilde T\into d\widetilde\theta
= 0,
\end{displaymath}
which means that $G_1$ preserves the contact form $\widetilde\theta$ and
therefore the Riemannian metric $g_{\tilde\theta}$.  Thus $G_1$ is
contained in the isometry group $\scr I(g_{\tilde\theta})$ of
$g_{\tilde\theta}$, which is compact.  Any sequence in $G_1$ therefore has
a subsequence that converges in the topology of $\scr I(g_{\tilde\theta})$,
which implies uniform convergence with all derivatives, and hence also in
the topology of $\scr A(M)$.  Since $G_1$ is closed in $\scr A(M)$, the
limit is again in $G_1$, so $G_1$ is compact, which is a contradiction.
Thus $f$ must vanish somewhere on $M$.

Now let $S$ denote the zero set of $f$.  In other words, $S$ is the set of
points of $M$ where $X$ is tangent to $H$.  At points of $S$, we have $X =
H^\theta_f$, so the assumption that $X$ does not vanish means that
$H^\theta_f\ne 0$ and therefore $df\ne 0$ along $S$.  Thus $S$ is a
compact, embedded hypersurface in $M$.  Moreover, since $H^\theta_f(f) =
df(H^\theta_f) = d\theta(H^\theta_f,H^\theta_f) = 0$, $X$ is tangent to
$S$.
\end{pf}

\begin{lemma}\label{onelemma}
The contact form $\theta$ can be chosen so that ${|\overline\partial_b
f|}_\theta = 1$ along $S$.  For any such $\theta$, the following relations
hold along $S$.
\begin{gather}
f_{;\alpha\beta} = 0;\label{rel1}\\
f_{;0} = 0;\label{rel2}\\
f^{;\beta} f_{;0\beta}  \hbox{ is real\rom{;}}\label{rel4}\\
f_{;0\alpha} = f^\beta f_{;0\beta} f_\alpha.\label{rel5}
\end{gather}
\end{lemma}

\begin{pf}
Since $H^\theta_f\ne 0$ at points of $S$, by rescaling $\theta$ we can
guarantee that the Levi norm ${|\overline\partial_b f|}_\theta$ is equal to
1 along $S$, as follows.  If $\widetilde\theta = u\theta$ for some positive
function $u\in C^\infty(M)$, then $\widetilde f = -\widetilde\theta(X) =
uf$, and so $\overline\partial_b \widetilde f = u \overline\partial_b f$
along $S$.  Therefore ${|\overline\partial_b {\widetilde
f}|}_{\tilde\theta}^2 = u^{-1}{|u \overline\partial_b f|}_{\theta}^2 =
u{|\overline\partial_b f|}_{\theta}^2$ along $S$, so replacing $\theta$ by
$u\theta$, where $u$ is any positive function that equals
${|\overline\partial_b f|}_\theta^{-2}$ along $S$, gives a contact form
satisfying the first statement of the lemma.

Let $\theta$ be any such contact form, and $f = -\theta(X)$.  Since $X$ is
a CR vector field, (\ref{fab}) evaluated along  $S$ gives (\ref{rel1}).

Next we prove (\ref{rel2}).  Since ${|\overline\partial_b f|}_\theta^2=1$
on $S$ and $X$ is tangent to $S$, at points of $S$ we have
\begin{align*}
0 &= X {|\overline\partial_b f|}_\theta^2
  = if^{;\bar\alpha}{(f_{;\bar\beta}f^{;\bar\beta})}_{;\bar\alpha} -
     if^{;\alpha}{(f_{;\bar\beta}f^{;\bar\beta})}_{;\alpha}\\
  &= if^{;\bar\alpha}f_{;\bar\beta}f^{;\bar\beta}{}_{\bar\alpha} -
     if^{;\alpha} f^{;\bar\beta}  f_{;\bar\beta\alpha}\\
  &= if^{;\bar\alpha} f^{;\beta} (f_{;\beta\bar\alpha} - f_{;\bar\alpha\beta})\\
  &= if^{;\bar\alpha} f^{;\beta} (i h_{\beta\bar\alpha} f_{;0})
  = -f_{;0}.
\end{align*}
This implies that $T$ is also tangent to $S$.  Therefore, on $S$ we have
\begin{displaymath}
0 = X T f = if^{;\bar\alpha}f_{;0\bar\alpha} - if^{;\alpha}f_{;0\alpha},
\end{displaymath}
which is (\ref{rel4}).

Finally, since $f_{;0}$ vanishes along $S$, we can write $f_{;0} = vf$ for
some smooth function $v$.  Differentiating, we find that
$f_{;0\alpha} = v f_{;\alpha}$ along $S$.  Then contracting with
$f^{;\alpha}$ yields $f_{;0\alpha}f^{;\alpha} = v f_{;\alpha}f^{;\alpha} =
v$, which proves (\ref{rel5}).
\end{pf}

\begin{lemma}\label{dthetalemma}
The contact form $\theta$ can be chosen so that $L_X\theta = L_Xd\theta=0$
at points of $S$.
\end{lemma}

\begin{pf}
Since $X$ is a contact vector field, $L_X\theta = -Tf\,\theta$ on all of
$M$.  Therefore, by (\ref{rel2}), $L_X\theta =0$ along $S$ provided
$\theta$ is chosen as in Lemma \ref{onelemma}.  If
$\widetilde\theta = v\theta$ is another contact form, then
\begin{align*}
 L_Xd\widetilde\theta &= d L_X\widetilde\theta \\
&= d( Xv\,\theta - v\,Tf\,\theta)\\
&= d(Xv)\wedge\theta + Xv\,d\theta - v\,d(Tf)\wedge\theta
-Tf\,dv\wedge\theta - v\,Tf\,d\theta.
\end{align*}
In order to preserve the property ${|\overline\partial_b f|}_\theta=1$
along $S$, we put $v = 1 + uf$ for some smooth function $u$ yet to be
determined.  Then we have the following relations at points of $S$:
\begin{gather*}
v=1;\\
Xv = 0;\\
d(Xv) = u\,d(Xf) + Xu\,df = u \,d(-fTf) + Xu\,df = Xu\,df.
\end{gather*}
(In the last line we have used the fact that $H_f^\theta f=0$.)
Therefore, along $S$ we have
\begin{align*}
 L_Xd\widetilde\theta &=  Xu\,df\wedge\theta -
d(Tf)\wedge \theta.
\end{align*}
To prove the lemma, it suffices to produce a real-valued function $u$ such
that $Xu\, df \equiv d(Tf)\pmod \theta$ along $S$.  In components this
means, using (\ref{rel5}),
\begin{displaymath}
(Xu) f_{;\alpha} = f_{;0\alpha} = f^{;\beta}f_{;0\beta}f_{;\alpha},
\end{displaymath}
which is equivalent to
\begin{displaymath}
Xu = f^{;\beta}f_{;0\beta}.
\end{displaymath}
The hypersurface $S$ is characteristic for this equation, so one might
expect that existence of a global solution would depend on the global
behavior of $f^{;\beta}f_{;0\beta}$ and the integral curves of $X$.
Surprisingly, however, it turns out that we can write down an explicit
solution.

To this end, consider the function $\Delta_bf= - (f_{;\beta}{}^{\beta} +
f_{;\bar\beta}{}^{\bar\beta})$.  Since $f_{;\beta}{}^{\beta} -
f_{;\bar\beta}{}^{\bar\beta} = nif_{;0}=0$ along $S$, we have the following
at points of $S$:
\begin{align}\label{xdeltabf}
X(\Delta_bf) &= - 2X(f_{;\beta}{}^{\beta})\\
&= - 2if^{;\bar\alpha} \tensor f{\down{;\beta}\up{\beta}\down{\bar\alpha}} +
     2if^{;\alpha}\tensor f{\down{;\beta}\up{\beta}\down{\alpha}}.\notag
\end{align}
Using the commutation relations (\ref{comm2}) and (\ref{comm3}), we can simplify the
third derivatives that appear above as follows:
\begin{align*}
\tensor f{\down{;\beta}\up{\beta}\down{\bar\alpha}}
&= {(f_{;\bar\beta}{}^{\bar\beta} + ni f_{;0})}_{;\bar\alpha}\\
&= f_{;\bar\beta\bar\alpha}{}^{\bar\beta}
   + i\delta_{\bar\alpha}{}^{\bar\beta}f_{;\bar\beta}{}_{0} -
   \tensor R {\down{\bar\beta}\up{\bar\rho}\down{ \bar\alpha} \up{\bar\beta}} f_{;\bar\rho}
   + n i f_{;0\bar\alpha} \\
&= {(iA_{\bar\beta\bar\alpha} f)}^{;\bar\beta} + if_{;\bar\alpha 0}
   - R_{\rho\bar\alpha}f^{;\rho} + ni f_{;0\bar\alpha}\\
&= iA_{\bar\beta\bar\alpha}f^{;\bar\beta} + if_{;\bar\alpha 0}
   - R_{\rho\bar\alpha}f^{;\rho} + ni f_{;0\bar\alpha}\\
&= (n+1)i f_{;0\bar\alpha} - R_{\rho\bar\alpha}f^{;\rho}.\\
\tensor f{\down{;\beta}\up{\beta}\down{\alpha}}
&= f_{;\beta\alpha} {}^{\beta} - i\delta_\alpha{}^\beta f_{;\beta 0}
   - \tensor R{\down\beta \up\rho \down\alpha \up\beta} f_{;\rho}\\
&= {(-iA_{\beta\alpha}f)}^{;\beta} - i f_{;\alpha 0}
   - R_{\alpha\bar\rho} f^{;\bar\rho}\\
&= -iA_{\beta\alpha}f^{;\beta} - i f_{;\alpha 0}
   - R_{\alpha\bar\rho} f^{;\bar\rho}\\
&= -if_{;0\alpha}
   - R_{\alpha\bar\rho} f^{;\bar\rho}.
\end{align*}
Inserting these relations into (\ref{xdeltabf}) and using (\ref{rel4}), we
obtain
\begin{displaymath}
X(\Delta_b f) =
   2(n+1)f^{;\bar\alpha} f_{;0\bar\alpha} + 2 f^{;\alpha} f_{;0\alpha}=
   2(n+2)f^{;\alpha} f_{;0\alpha}.
\end{displaymath}
Thus the conclusion of the lemma holds if we replace $\theta$ by
$\widetilde\theta = (1+uf)\theta$, where $u$ is any function that is equal
to $\Delta_b f/(2(n+2))$ along $S$.
\end{pf}

\begin{pf*}{Proof of Theorem \ref{fixptthm}}
Suppose $G_1$ has no fixed points, and let $\{\phi_j\}\subset G_1$ be any
sequence; we will show it has a convergent subsequence, which is a
contradiction.

If $\theta$ is chosen as in Lemma \ref{dthetalemma}, the Webster metric
$g_\theta = d\theta(\cdot, J_\theta\cdot) + \theta^2$ is preserved along
$S$ by $G_1$, since $L_X\theta = L_X d\theta = L_XT = L_X J_\theta=0$ along
$S$.  Since the group of isometries of ${\left.g_\theta\right|}_S$ is
compact, there is a subsequence, still denoted $\{\phi_j\}$, whose
restrictions to $S$ converge uniformly with all derivatives.

I claim that the two-jets of the sequence $\{\phi_j\}$ converge at all
points of $S$.  This implies that $\{\phi_j\}$ converges in the topology of
$\scr A(M)$ by the following standard argument.  Let $Y\to M$ be the Chern
CR structure bundle of $M$ \cite{CM}; then every CR automorphism $\phi$ of
$M$ lifts naturally to an automorphism $\widetilde \phi$ of $Y$ preserving
the Chern connection.  Moreover, choosing any point $\xi\in Y$ over $S$,
the orbit map $\phi\mapsto\widetilde \phi(\xi)$ gives a closed embedding of
$G_1$ into $Y$ \cite[Thm.\ I.3.2]{Kobayashi}.  Convergence of the two-jets
of $\{\phi_j\}$ along $S$ implies that $\{\widetilde
\phi_j(\xi)\}$ converges in $Y$, which therefore implies that $\{\phi_j\}$
converges to some element $\phi\in G_1$.

To see that convergence of $\phi_j$ on $S$ implies convergence on the
two-jet level, let $Y$ denote the vector field $J_\theta X = J H_f^\theta$.
Since $Y = f^{;\bar\alpha}Z_{\bar\alpha} + f^{;\alpha} Z_{\alpha}$ by
(\ref{hf}), it follows that $Yf = 2 f^{;\bar\alpha} f_{;\bar\alpha} = 2$
along $S$.  In particular, $Y$ is always transverse to $S$.  For any
element $\phi\in G_1$, using the facts that $\phi$ is a CR automorphism and
$\phi_*X=X$, we compute
\begin{equation}\label{phistary}
\begin{aligned}
\phi_*Y &= \phi_*(J H_f^\theta) = J (\phi_* H_f^\theta)\\
&= J( \phi_*(X + f T) )
= J_\theta X + (f\circ \phi^{-1}) J_\theta (\phi_*T )\\
&= Y + (f\circ\phi^{-1}) J_\theta (\phi_*T).
\end{aligned}
\end{equation}

Given a point $p\in S$, we can choose coordinates $(x_1,\dots,x_{2n},y)$ on
some neighborhood $U$ of $p$ such that $y=0$ on $S$ and $Y
\equiv\partial/\partial y$ .  In these coordinates, $f(x,y) = 2 y$, since
both sides satisfy the ordinary differential equation $Yf=2$ with initial
condition $f=0$ on $S$.  Choosing analogous coordinates $(\widetilde
x,\widetilde y)$ near $\lim_{j\to \infty}\phi_j(p)\in S$ and taking $j$
sufficiently large and $U$ sufficiently small, we may assume that
$\phi_j(x,0)$ lies in a single coordinate chart for all $j$ and all
$(x,0)\in U\cap S$.  In these coordinates, we can consider $\phi_j$ as a
vector-valued function of $(x,y)$.  We already know that all the
$x$-derivatives of $\phi_j$ converge on $S$, so we need only consider the
$y$-derivatives.

In coordinates, (\ref{phistary}) becomes
\begin{displaymath}
\frac {\partial\phi_j(x,y)} {\partial y} =
(0,1) + 2 y V_j(x,y),
\end{displaymath}
where $V_j(x,y)$ is the coordinate representation of $J_\theta
((\phi_j)_*T_{(x,y)})$.  In particular, along $S$ we get $\partial
\phi_j/\partial y\equiv(0,1)$,
which certainly converges together with all its $x$-derivatives.  Moreover,
the values of ${\partial^2\phi_j}/ {\partial y^2} = 2V_j = J_\theta
((\phi_j)_*T)$ along $S$ are determined by the 1-jets of $\phi_j$ along
$S$, and we have already shown these converge on $S$.  Therefore
all second derivatives of $\phi_j$ converge  on $S$.  This
completes the proof of Theorem
\ref{fixptthm}.
\end{pf*}
\vfill

\newpage\raggedbottom

\bigskip
\noindent\fbox{\sc Version 2.2 -- May 2, 1994}

\end{document}